\documentclass[12pt]{amsart}
 \title[Representation categories of matrix quantum groups]{On the representation categories of matrix quantum groups of type $A$}
\author[P.H.Hai]{Ph\`ung  H$\grave{\mbox{\^o}}$  Hai}
\address{ Hanoi Institute of Mathematics \\ P.O. Box 10000, Bo Ho, Hanoi}
 \email{ phung@math.ac.vn}
\curraddr{Department of Mathematics, University of Duisburg-Essen, 45117 Essen, Germany}
\subjclass[2000]{Primary 16W30, 17B37, Secondary 17A45, 17A70}
\thanks{This work was supported in part by the Nat. Program for Basic 
Sciences Research of Vietnam and the   ``DFG-Schwerpunkt
 Komplexe Mannigfaltigkeiten''.}

\newtheorem{thm}{Theorem}[section]
\newtheorem{lem}[thm]{Lemma}

\newtheorem{cor}[thm]{Corollary}

\def\proof{{\it Proof.\ }}
\def\H{{\mathcal H}}
\newcommand{\Ss}{\mathfrak S}

\def\id{{\mathchoice{\mbox{\rm id}}
                    {\mbox{\rm id}}
                    {\mbox{\scriptsize\rm id}}
                    {\mbox{\tiny\rm id}} }}

\def\op{{\mathchoice{\mbox{\rm op}}
                    {\mbox{\rm op}}
                    {\mbox{\scriptsize\rm op}}
                    {\mbox{\tiny\rm op}} }}

\def\exterior{{\mathchoice{\mbox{\large$\wedge$}}
{\mbox{\large$ \wedge$}}
{\mbox{\small$\wedge$}}
{\mbox{\scriptsize$\wedge$}} }}
\def\sym{{\mathchoice{\mbox{\sf S\hspace{0.1ex}}}
{\mbox{\sf S\hspace{0.1ex}}}
{\mbox{\scriptsize\sf S\hspace{0.1ex}}}
{\mbox{\tiny\sf S\hspace{0.1ex}}} }}

\begin{document}
\begin{abstract} A quantum groups of type $A$ is defined in terms of
a Hecke symmetry. We show in this paper that the representation
category of such a quantum group  is uniquely determined as an abelian braided monoidal
category by the bi-rank of the Hecke symmetry.\end{abstract}
\maketitle
\bibliographystyle{plain}

\section{Introduction}
A matrix quantum group of type $A$ is defined as the ``spectrum'' of
the Hopf algebra associated to a closed solution of the (quantized)
Yang-Baxter equation and the Hecke equation (called a Hecke
symmetry). Explicitly, let $V$ be a  vector space
(over a field) of finite dimension $d$. An invertible operator 
$R:V \otimes V \longrightarrow
V \otimes V$ is called a Hecke symmetry if it satisfies  the equations
\begin{equation}\label{eq01}R_1R_2R_1=R_2R_1R_2, \quad \mbox{where
} R_1:=R \otimes \id_V, R_2:=\id_V \otimes R \end{equation}
(the Yang-Baxter equation),
\begin{equation}\label{eq02} (R+1)(R-q)=0, \quad q \neq 0;-1\end{equation}
(the Hecke equation) and is closed in the sense that the half dual operator
$$ R^\sharp:V^*\otimes V \longrightarrow V \otimes V^*;\quad \langle
R^\sharp(\xi \otimes v),w \rangle=\langle \xi,R(v \otimes w)\rangle$$
is invertible.

Given such a Hecke symmetry one constructs a Hopf algebra $H$
as follows. Fix a basis $\{x_i;1\leq i \leq d \}$ of $V$ and let $R^{i j}_{k l}$ 
be the matrix of $R$ with respect to this basis. 
As an algebra $H$ is generated by two sets of generators
$\{z^i_j,t^i_j;1\leq i \leq d \}$, subject to the following relations
(we will always adopt the convention of summing over the indices
that appear in both upper and lower places):
$$\begin{array}{l}
R^{i j}_{p q}z^p_k z^q_l=z^{i}_{m}z^j_n R^{m n}_{kl}\\
z^i_k t^k_j=t^i_k z^k_j=\delta^i_j\end{array}$$
In case $R$ is the usual symmetry operator: $R(v\otimes w)=w\otimes v$
(thus $q=1$), $H$ is isomorphic to the function algebra on the algebraic
group $GL(V)$. 

The most well-known Hecke symmetry is the Drinfeld-Jimbo
solutions of series $A$ to the Yang-Baxter equation 
(fix a square root  $\sqrt q$ of $q$)
\begin{equation}\label{eq03}
R^{d}_q(x_i \otimes x_j)=\left[\begin{array}{ll}
q x_i \otimes x_i & \mbox{if } i=j \\
\sqrt q x_j \otimes x_i & \mbox{if }i >j \\
 \sqrt q x_j \otimes x_i- (q-1) x_i \otimes x_j & \mbox{if } i<j \end{array}\right.
\end{equation}
In the ``classical" limit $q \rightarrow 1$, $R^{d}_q$ reduces to the
usual symmetry operator. There is also a super version of these 
solutions due to Yu.~Manin \cite{manin2}. Let $V$ be a  vector superspace of 
super-dimension $(r|s)$, $ r+s=d$, and
let $\{x_i\}$ be a homogeneous basis of $V$, the parity of $x_i$ is denoted
by $\hat i$. The Hecke symmetry $R_q^{r|s}$ is given by
\begin{equation}\label{eq04}
R^{r|s}_q(x_i \otimes x_j)=\left[\begin{array}{ll}
(-1)^{\hat i}q x_i \otimes x_i & \mbox{if } i=j \\
(-1)^{\hat i\hat j}\sqrt q x_j \otimes x_i & \mbox{if }i >j \\
 (-1)^{\hat i\hat j}\sqrt q x_j \otimes x_i- 
 (q-1) x_i \otimes x_j & \mbox{if } i<j \end{array}\right.
\end{equation}
In the ``classical" limit $q \rightarrow 1$, $R^{r|s}_q$ reduces to the
super-symmetry operator $R(x_i\otimes x_j)=(-1)^{\hat i \hat j}x_j \otimes x_i.$

The quantum group associated to the Drinfeld-Jimbo solution (\ref{eq03})
 is called the standard quantum deformation of the general linear group
 or simply standard quantum general linear group. Similarly, the
 quantum general linear super-group is determined in terms of the solution
 (\ref{eq04}) (actually, some signs must be inserted in the definition, see
 \cite{manin2} for details).
 
 There are many other non-standard Hecke symmetries and there is so far
 no classification of these solutions except for the case the dimension of
 $V$ is 2. On the other hand, many properties of the associated quantum
 groups to these solutions are  obtained in an abstract way. 
 The aim of this work is to study representation category of the matrix quantum
 group associated to a Hecke symmetry, by this we understand the comodule
 category over the corresponding Hopf algebra.
 The pair $(r,s)$, where $r$ is the number of roots and $s$ is the number of
 poles of $P_\exterior(t)$ (see \ref{sect1.2.4}), is called the {\em bi-rank} of the Hecke symmetry.
  The main result of this paper is that the category of comodules over the Hopf 
 algebra associated to a Hecke symmetry as a braided monoidal abelian category
 depends only on the parameter $q$ and the bi-rank.
 
 The proof of the main result is inspired by the work \cite{bichon} of J.~Bichon, 
 whose idea
 was to use a result of P.~Schauenburg on the relationship between equivalences
 of comodule categories a pair of Hopf algebras and bi-Galois extensions.
 
 The main result implies that the study of representations of a matrix 
 quantum group of type $A$ can be reduced to the study of that of a standard quantum
general linear group. The latter has been studied by R.~Zhang \cite{zhang}.
In particular we show that the homological determinant is always one-dimensional.

This work was carried out during the author's visit at the Department
of Mathematics, University of Duisburg-Essen. He would like to thanks Professors
H. Esnault and E. Viehweg for the financial support through their
Leibniz-Preis  and for their hospitality.

\section{Matrix quantum group of type $A$}\label{sect1}
Let $V$ be a vector space  of finite dimension $d$ over a field $k$ of characteristic zero. 
Let $R:V \otimes V \longrightarrow V \otimes V$ be  a {\it Hecke symmetry}.
Throughout this work we will  assume that $q$ {\em is not a root of
 unity other then the unity itself}. The matrix ${R^\sharp}$ 
 is given by ${R^\sharp}^{kl}_{ij}=R_{jl}^{ik}$.
Therefore, the invertibility of $R^\sharp$ can be expressed as follows: 
there exists a matrix $P$ such that  $P^{im}_{jn}R^{nk}_{ml}=\delta^i_l\delta^k_j$.
Define the following algebras:
\begin{eqnarray*} 
 \sym&:=&k \langle x_1,x_2,  \ldots,x_d \rangle/(x_kx_lR^{kl}_{ij}=qx_ix_j) \\
 \exterior&:=&k \langle x_1,x_2,  \ldots,x_d \rangle/(x_kx_lR^{kl}_{ij}=-x_ix_j) \\
 E&:=&k \langle z^1_1,z^1_2,  \ldots ,z^d_d \rangle/(z^{i}_{m}z^j_nR^{mn}_{kl}=
R^{ij}_{pq}z^p_kz^q_l) \\
H&:=& k \langle z^1_1,z^1_2, \ldots ,z^d_d,t^1_1,t^1_2,\ldots,t^d_d\rangle
\raisebox{-1ex}{$\left/\left(\begin{array}{l}
z^{i}_{m}z^j_n R^{m n}_{kl}=R^{i j}_{p q}z^p_k z^q_l\\
z^i_k t^k_j=t^i_k z^k_j=\delta^i_j\end{array}\right)\right.$}
\end{eqnarray*}
where $\{x_i \}$, $\{z^i_j \}$ and $\{t^i_j \}$ are sets of generators.

The algebras $\exterior$ and $\sym$ are called the {\em quantum anti-symmetric} and 
{\em quantum symmetric algebras} associated to $R$. Together they define a
quantum vector space.

 The algebra $E$ is in fact  a bialgebra with coproduct and counit given by
$$\Delta(z^i_j)=z^i_k \otimes z^k_j, \quad \varepsilon(z^i_j)=\delta^i_j.$$ 
The algebra $H$ is a Hopf
algebra with $\Delta(z^i_j)=z^i_k \otimes z^k_j$, $\Delta(t^i_j)=z_j^k \otimes z_k^i$, 
$\varepsilon(z^i_j)=\varepsilon(t^i_j)=\delta^i_j$. For the
antipode, let $C^i_j:=P^{i m}_{j m}$. Then $S(z^i_j)=t^i_j$ and
\begin{equation}\label{eq05}S(t^i_j)=C^i_k z^k_l {C^{-1}}^l_j
\end{equation}
 \cite[Thm. 2.1.1]{ph97b}. The matrix $C$ plays an important role
 in our study, its trace is called the quantum rank of the Hecke symmetry
$\mbox{Rank}_qR:=\mbox{tr }(C)$, see \ref{sect1.3.2}.

The bialgebra $E$ is considered as {\em the function 
algebra on a quantum semi-group of type
$A$} and the Hopf algebra $H$ is considered  as {\em the function algebra on a 
matrix quantum groups of $A$}. Representations of this 
(semi-)group are thus comodules over $H$ (resp. $E$).

\subsection{Comodules over $E$}\label{sect1.2}
The space  $V$ is a comodule over $E$
by the map $\delta:V \longrightarrow V \otimes E; x_i \longmapsto x_j \otimes z^j_i$. 
Since $E$ is a bialgebra, any  tensor power of $V$
is also a comodule over $E$. The
map $R:V \otimes V \longrightarrow V \otimes V$ is a  comodule map. 
The classification of $E$-comodules is 
done with the help of the action of the Hecke algebra.

\subsubsection{The Hecke algebras} \label{sect1.2.1}
The Hecke algebra $\H_n=\H_{q,n}$ 
has generators $t_i, 1\leq i\leq n-1$, subject to the relations: 
\begin{eqnarray*}&&t_it_j=t_jt_i,  |i-j|\geq 2;\\ 
&& t_it_{i+1}t_i=t_{i+1}t_it_{i+1}, 1\leq i\leq n-2;\\
&& t_i^2=(q-1)t_i+q.\end{eqnarray*} 
There is a $k$-basis in $\H_n$ indexed by permutations of $n$ elements
: $t_w,w \in \Ss_n$ ($\Ss_n$ is the permutation group),
 in such a way that $t_{(i,i+1)}=t_i$ and
$t_w t_v=t_{w v}$ if the length of $w v$ is equal to the sum of the 
length of $w$ and the length of $v$. 

If $q$ is not a root of unity of
degree greater than 1, $\H_n$ is a semisimple algebra.
It is isomorphic to the direct product of its minimal two-sided ideals,
which are themselves simple algebras. The minimal two-sided ideals
can be indexed by partitions of $n$. Thus
$$\H_n \cong \prod_{\lambda \vdash n}A_\lambda$$
$A_\lambda$ denotes the minimal two-sided ideal corresponding to $\lambda$.
Each $A_\lambda$ is a matrix ring over the ground field and one
 can choose a basis $\{e_\lambda^{i j}; 1\leq i,j \leq d_\lambda \}$
 such that 
$$e_\lambda^{i j}e_\lambda^{k l}=\delta^j_k e^{il}_\lambda$$
$d_\lambda$ is the dimension of  the simple $\H_n$-comodule corresponding
to $\lambda$ and can be computed by the combinatorics of
$\lambda$-tableaux. In particular, $\{ e^{i i}_\lambda, 1\leq i\leq d_\lambda \}$
 are mutually orthogonal conjugate primitive idempotents of $\H_n$.
For more details, the reader is referred to
\cite{dj1,dj2}.

\subsubsection{An action of $\H_n$}\label{sect1.2.2}The Hecke symmetry
$R$ induces an action of the Hecke algebra $\H_n=\H_{q,n}$ on
$V^{\otimes n}$, $t_i \longmapsto R_i=\id^{i-1}\otimes R \otimes\id^{n-i-1}$ 
which commutes with the coaction of $E$. 
The action of $t_w$ will be denoted by $R_w$.

 Thus, each element of $\H_n$ determines an endomorphism of $V^{\otimes n}$
as $E$-comodule. For $q$ not a root of unity of degree greater 1, 
the converse is also true: each endomorphism
of $V^{\otimes n}$ represents the action of an element of $\H_n$, moreover
 $V^{\otimes n}$ is semi-simple and its simple subcomodules
can be  given as the images of the endomorphisms determined by
primitive idempotents of $\H_n$, conjugate idempotents
(i.e. belonging to the same minimal two-sided ideal) determine
isomorphic comodules \cite{ph97b}. 

Since conjugate classes of primitive idempotents
of $\H_n$ are indexed by partitions of $n$, simple subcomodules of 
$V^{\otimes n}$ are  indexed by a subset of partitions of $n$. Thus $E$ is
cosemisimple and its simple comodules are indexed by a subset of
partitions. 

\subsubsection{ Quantum symmetrizers}\label{sect1.2.3}
Denote $$[n]_q:=\displaystyle\frac{q^n-1}{q-1}$$
 The primitive idempotent 
corresponding to partition $(n)$ of $n$,
$$X_n:=\displaystyle\frac1{[n]_{q}}\sum_{w\in \Ss_n}R_w$$
 determines a simple
comodule isomorphic to the $n$-th homogeneous component $\sym_n$
of the quantum symmetric algebra $\sym$  (the $n$-th quantum symmetric power)
and the primitive idempotent corresponding to partition $(1^n)$ of $n$,
$$Y_n:=\displaystyle\frac1{[n]_{1/q}}\sum_{w\in \Ss_n}(-q)^{-l(w)}R_w$$ 
determines a simple comodule isomorphic to the $n$-th homogeneous
component $\exterior_n$ of the quantum exterior algebra $\exterior$
(the $n$-th quantum anti-symmetric power).

\subsubsection{The bi-rank}\label{sect1.2.4}
There is a determinantal formula in the Grothendieck
ring of finite dimensional $E$-comodules which computes simple comodules
in terms of quantum symmetric tensor powers \cite{ph97b}:
\begin{equation}\label{eq06}
I_\lambda=\det |\sym_{\lambda_i-i+j}|_{1\leq i, j \leq k};
\quad k\mbox{ is the length of }\lambda.
\end{equation}
Consequently we have a similar form for the dimensions of simple comodules.
It follows from this and a theorem of Edrei on P\'olya frequency sequences
that the Poincar\'e series of  $\exterior$ is rational function with
negative zeros and positive poles \cite{ph97c}. The pair $(r,s)$ where
$r$ is the number of zeros and $s$ is the number of poles is called
the {\em bi-rank} of the Hecke symmetry.
 It then follows from (\ref{eq06}) that
the $E$-comodule $I_\lambda$ is non-zero, and hence simple,
if and only if $\lambda_r \leq s$, the set of partitions of $n$ satisfying this property is
denoted by $\Gamma^{r,s}_n$. Simple $E$-comodules are thus completely
classified in terms of the bi-rank. 

\subsection{The Hopf algebra $H$ and its comodules}
\subsubsection{The Koszul complex}\label{sect1.3.2}
Through the natural map $E \longrightarrow H$ 
 $E$-comodules are comodule over $H$. Since $H$ is a Hopf algebra, for
 the comodules $\sym_n$ and $\exterior_n$,  
 their dual $\sym_{n}{}^*,\exterior_n{}^* $ space are also comodules over $H$.
One can define $H$-comodule maps
 $$d^{k,l}: K^{k,l}:= \exterior_k \otimes \sym_{l}{}^*{\longrightarrow}
K^{k+1,l+1}:=\exterior_{k+1}\otimes \sym_{l+1}{}^*,\quad k,l\in\mathbb Z,$$
 in such a way that the sequence
$$K^{a}: \cdots {\longrightarrow}\exterior_{k-1}\otimes \sym_{l-1}{}^*
{\longrightarrow}\exterior_k \otimes \sym_{l}{}^*{\longrightarrow}
\exterior_{k+1}\otimes \sym_{l+1}{}^*\cdots$$
$(a=k-l)$ is a complex. This complexes were introduced by Manin for the case of
standard Hecke symmetry \cite{manin2} and studied by 
Gurevich, Lyubashenko, Sudbery \cite{gur1,ls}. 

It is expected that the homology of this complexes is concentrated at a certain
term where it has dimension one, in this case it induce a group-like
element in $H,$ called {\em the homological
determinant} as suggested by Manin.  Gurevich showed that all the complexes $K^a$, 
$a \in \mathbb Z$, except might be for the complex $K^b$ with 
$[-b]_q=-\mbox{rank}_q R$, are exact.  For the case of even Hecke symmetries
he showed that the homology is one-dimensional~\cite{gur1}. 
The homology of the complex $K^b$
was shown to be one-dimensional by Manin for the case of standard Hecke
symmetry~\cite{manin2}. This fact has also been shown
by Lyubashenko-Sudbery for Hecke sums of odd and even Hecke symmetries~\cite{ls}.
A combinatorial proof for Hecke symmetries of birank $(2,1)$ was given in
\cite{ph02d}.

In \cite{ph02c} the author showed that the homology should be non-vanishing
at the term $\exterior_r \otimes \sym_{s}{}^*$ and consequently the quantum rank
$\mbox{rank}_q R=\mbox{tr }(C)$ is equal to $-[s-r]_q$.

\subsubsection{The integral}\label{sect1.3.1}
 In the study of the category $H$-comod, the integral over $H$ 
 plays an important role as shown in \cite{ph02d}. By definition, a
right integral over $H$ is a (right) comodule map $H \longrightarrow k$
where $H$ coacts on itself by the coproduct and on the base field $k$
by the unit map. 
The existence of the integral on $H$ was proven in \cite[Thm.3.2]{ph98b},
under the assumption that
$\mbox{rank}_qR=-[s-r]_q$, which
was later shown in \cite{ph02c} for an arbitrary Hecke symmetry.
  In fact, an explicit form for the integral was given. Since
we will need it later on, let us recall it here.

For a partition $\lambda$ of $n$, let $[\lambda]$ be the corresponding tableau
 and for any node $x \in[ \lambda]$, $c(x)$ be its contents,
$h(x)$ its hook-length,
 $n(\lambda):=\sum_{x \in[\lambda]}c(x)$ (see \cite{mcd2} for details). 
 
 Let 
\begin{eqnarray*}
 p_\lambda := \prod_{x \in [\lambda]\setminus[(s^r)]}
q^{r-s}{[c(x)+r-s]_q}^{-1},\quad
k_\lambda := q^{n(\lambda)}\prod_{x \in[\lambda]}[h(x)]_q^{-1},
\end{eqnarray*}
where $(s^r)$ is the sub-tableau of $\lambda$ consisting of nodes
in the $i$-th row and $j$-th column with $i \leq s, j \leq r$. In particular,
$p_\lambda=0$ if $\lambda_r<s$. Let $\Omega^{r,s}_n$ denote the set of partitions
 from $\Gamma^{r,s}_n$ such that $p_\lambda \neq 0$. 
Thus $\Omega^{r,s}_n=\{\lambda \vdash n; \lambda_r=s\}$.

Denote for each set of indices $I=(i_1,i_2, \ldots, i_n)$, $J=(j_1,j_2,\ldots,j_n)$
$$ Z^I_J:=z^{i_1}_{j_1}z^{i_2}_{j_2}\ldots z^{i_n}_{j_n};\quad 
 T^{I'}_{J'}:=t^{i_n}_{j_n} \ldots t^{i_2}_{j_2}t^{i_1}_{j_1} .$$
 Then the value of the integral on $Z^I_JT^{K'}_{L'}$ can be given as follows
 \begin{equation}\label{eq07}
 \int(Z_I^J T_{K'}^{L'})=
\sum_{\lambda \in \Omega^{r,s}_n \atop {1\leq i,j\leq d_\lambda}}
\frac{p_\lambda}{k_\lambda}(C^{\otimes n}E_\lambda^{ij})_I^L (E_\lambda^{ji})_K^J
 \end{equation}
where $E^{i j}_\lambda$ is the matrix of the basis element $e^{i j}_\lambda$ in
the representation $\rho_n$, the matrix $C$ is given in (\ref{eq05}).
 In particular, the left hand-side is zero if $n < r s$.

\section{Bi-Galois extensions}
Let $A$ be a Hopf algebra over a field $k$. A right $A$-comodule algebra 
is a right $A$-comodule with the structure of an algebra on it such that the
structure maps (the multiplication and the unit map) are $A$-comodule maps.
A right $A$-Galois extension $M/k$ is a right $A$-comodule algebra $M$
such that  the Galois map 
\begin{equation}\label{eq08}
\kappa_r:M \otimes M \longrightarrow M \otimes A;
\quad \kappa_r(m \otimes n)=\sum_{(n)}m n_{(0)}\otimes n_{(1)}\end{equation}
is bijective. Similarly one has the notion of left $A$-Galois extension, in which $M$ is
a left $A$-comodule algebra and the Galois map is 
$\kappa_l:M \otimes M \longrightarrow A \otimes M;$  $m \otimes n \longmapsto
\sum_{(m)}m_{(-1)} \otimes m_{(0)}n$.

\begin{lem}\label{lem1} Let $M$ be a right $A$-comodule algebra. Assume
that there exists an algebra map $\gamma:A \longrightarrow M^{\rm op} \otimes M,$
$a \longmapsto \sum_{(a)}a^- \otimes a^+$ such that
 the following equations in $M \otimes M$ hold true
\begin{equation}\label{eq08a}\begin{array}{l}
  \sum_{(m)}m_{(0)} m_{(1)}{}^- \otimes m_{(1)}{}^+= 1 \otimes m, \quad
m \in M\\
 \sum_{(a)} a^- a^+{}_{(0)} \otimes a^+{}_{(1)}= 1 \otimes a;\quad
a \in A \end{array}\end{equation}
Then $M$ is a right $A$-Galois extension of $k$.
For left Galois extension the conditions read: 
$\gamma: A \longrightarrow M \otimes M^{\op}$, 
$a \longmapsto \sum_{(a)}a^+ \otimes a^-,$
\begin{equation}\label{eq08b}\begin{array}{l}
 \sum_{(m)}m_{(-1)}{}^+ \otimes m_{(-1)}{}^- m_{(0)}=m \otimes 1;\quad m \in M \\
 \sum_{(a)}a^+{}_{(-1)} \otimes a^+{}_{(0)} a^-= a \otimes 1; \quad a \in A
\end{array}\end{equation}
\end{lem}
\proof The inverse to $\kappa_r$ is given in terms of $\gamma$ as follows:
$m\otimes a\longmapsto \sum_{(a)} m a^-\otimes a^+$. For
$\kappa_l$, the ivnerse is given by 
$a\otimes m\longmapsto \sum_{(a)} a^+\otimes a^- m$. $\Box$

\bigskip
\noindent {\sc Remark.} We see from the proof 
that the map $\gamma$ can be obtained
from $\kappa_r$ as follows:
$\gamma(a)=\kappa_r{}^{-1}(1\otimes a)$. Then, one can show
that $\gamma$ is an algebra homomorphism.
In fact, in the above proof, we do not use the fact that
$\gamma$ is an algebra homomorphism. We assume it
however, since the equations in (\ref{eq08a})
and (\ref{eq08b}) respect the multiplications in
$A$ and $M$, that is, if an equation holds true for $a$ and $a'$ in $A$
or $m$ and $m'$ in $M$ then it holds true for the products $a a'$ or $m m'$
respectively. Therefore it is sufficient  to check this conditions on a set 
of generators of $A$ and $M$.
\bigskip

Let now $A$ and $B$ be Hopf algebras and $M$ an $A-B$-bi-comodule, i.e.
$M$ is left $A$-comodule and right $B$-comodule and the two coactions
are compatible. $M$ is said to be an $A-B$-bi-Galois extension of $k$ if it is
both a left $A$-Galois extension and a right $B$-Galois extension of $k$. We
will make use of the following fact \cite[Cor.5.7]{schauen3}:

{\em There exists a 1-1 correspondence between the set of isomorphic classes
of $A-B$-bi-Galois extension of $k$ and $k$-linear monoidal equivalences 
between the categories of comodules over $A$ and $B.$}

The equivalence functor is given in terms of the co-tensor product with the
bi-comodule. Recall that each $A-B$-bi-comodule $M$ defines an additive functor
from the category $A$-comod of right $A$-comodules to the category $B$-comod
$X \longmapsto X \Box_A M$, where the co-tensor product $X \Box_A M$ is defined
as the equalizer of the two maps induced from the coactions on $A$:
$$X \Box_A M \longrightarrow X \otimes_k M 
\raisebox{.5ex}{$\begin{array}{c}{{}_{\id \otimes \delta_M}}\\[-1ex]
\longrightarrow \\[-1.5ex] 
\longrightarrow \\[-2ex]
{}_{\delta_X \otimes \id}  \end{array}$} 
X \otimes_k A \otimes_k M,$$
or, explicitly, 
$$X\Box_AM=\{ x\otimes m\in X\otimes_kM| 
\sum_{(m)}x\otimes m_{(-1)}\otimes m_{(0)}=\sum_{(x)}x_{(0)}\otimes 
x_{(1)}\otimes m\}.$$
The coaction of $B$ on $X \Box_A M$ is induced from that on $M.$

\section{A bi-Galois extension for matrix quantum groups}
Let $R$ and $\bar R$ be Hecke symmetries and $H$, $\bar H$
be the associated Hopf algebras. We construct in this subsection
an $H-\bar H$-bi-Galois extension.

Assume that $R$ is defined on a vector space of dimension $d$ and
$\bar R$ is defined over a vector space of dimension $\bar d$. Consider
the algebra $M=M_{R,\bar R}$ generated by elements 
$a^i_\lambda, b^\lambda_i; 1 \leq i \leq d, 1 \leq \lambda \leq \bar d$, 
subject to the following relations
 \begin{eqnarray*} 
&&R^{i j}_{p q}a^p_\lambda a^q_\mu =
a^{i}_{\nu}a^j_\gamma \bar R^{\nu \gamma}_{\lambda \mu},\\
 && a^i_\lambda b^\lambda_j=\delta^i_j;\quad 
 b^\lambda _k a^k_\mu= \delta^\lambda_\mu  
\end{eqnarray*}
The following equations can also be deduced from the equations above
\begin{eqnarray*}
&& R^{m n}_{kl}b^{\lambda}_{n}b^\mu_m = 
b^\gamma_k b^\nu_l\bar R^{\lambda \mu}_{\nu \gamma}\\
&& P^{qp}_{l k}a_\nu^l b^\gamma_q =  
b_k^\mu a^p_\lambda \bar P^{\gamma \lambda}_{\nu \mu}\\
&& a_\gamma^l C^q_l b^\nu_q  =  \bar C^\nu_\gamma 
\end{eqnarray*}
The proof is completely similar to that of \cite[Thm.2.1.1]{ph97b}.

\begin{lem}\label{lem2} Assume that the algebra $M$ constructed above is
non-zero. Then it is an $H-\bar H$-bi-Galois extension of $k$.\end{lem}
\proof  The coactions of $H$ and $\bar H$ on $M$ are given by
\begin{eqnarray*} \delta &:& M \longrightarrow H \otimes M; \quad
   a^i_j \longmapsto z^i_k \otimes a^k_j, \quad
  b^j_i \longmapsto t^k_i \otimes b^j_k  \\
 \bar \delta &:& M \longrightarrow M \otimes \bar H; \quad
  a^i_j \longmapsto a^i_k \otimes \bar z{}^k_j, \quad
  b^j_i \longmapsto  b_i^k \otimes \bar t{}^j_k
 \end{eqnarray*}

 The verification that this maps induce a structure of left $H$-comodule
 (resp. right $\bar H$-comodule) algebra over $H$ and an $H- \bar H$ 
 bi-comodule structure is
 straightforward. 
 
According to Lemma \ref{lem1} and the remark following it, to show
that $M$ is a left $H$-Galois extension of $k$
it suffices to construct the map $\gamma$ satisfying the condition of the lemma. Define
\begin{eqnarray*}
&& \gamma(z^i_j)= a^i_\mu \otimes b^\mu_j\\
&& \gamma(t^i_j)= b_j^\mu \otimes \bar C^{-1}{}^\nu_\mu a^l_\nu C^j_l
\end{eqnarray*}
and extend them to algebra maps.
 Using the relations on $M$ one can check easily that this map gives rise to
 an algebra homomorphism $H \longrightarrow M \otimes M^{\rm op}$.
 Since $M$ is now an $H$-comodule algebra and since $\gamma$ is algebra
 homomorphism, the equations in Lemma \ref{lem3} respect the multiplications in $M$
 and in $H$, that is, it suffices to check them for the generators $z^i_j$ and
 $t^i_j$ which follows immediately from the relations mentioned above on
 the $a^i_\lambda$ and $b^\mu_j$. \hfill  $\Box$
 
 Notice that in the proof of this lemma the Hecke equation is not used.
 
 \begin{lem}\label{lem3}
Let $R$ and $\bar R$ be Hecke symmetries defined over 
 $V$ and $\bar V$ respectively. Assume that they are defined for
 the same value $q$ and have the same bi-rank. Then the associated
 algebra $M=M_{R,\bar R}$ is non-zero.\end{lem}
 \proof To show that $M$ is non-zero we construct a linear functional
 on $M$ and show that this linear functional attains a non-zero value
 at some element of $M.$ The construction of the linear functional
 resembles the integral on the Hopf algebra $H$ given in the previous 
 section. In fact, using the same method as in the proof of Theorem 3.2 and
 Equation 3.6 of \cite{ph98b} we can show that there is a linear functional
 on $M$ given by
   \begin{eqnarray*}
 \int(A_\Lambda^J B_{K'}^{\Gamma'})=
\sum_{\lambda \in \Omega^{r,s}_n \atop {1\leq i,j\leq d_\lambda}}
\frac{p_\lambda}{k_\lambda}(\bar C^{\otimes n}
\bar E_\lambda^{i j})_\Lambda^\Gamma (E_\lambda^{ji})_K^J
 \end{eqnarray*}
 where $\Lambda,\Gamma,I,J$ are multi-indices of length $n$ and $(r,s)$
 is the bi-rank of $R$ and $\bar R$.

According to  \ref{sect1.2.4} and \ref{sect1.3.1}
 for $n \geq rs$ and $\lambda \in \Omega^{r,s}_n$ the
matrices $E_\lambda^{ji}$ and $\bar E_\lambda^{ij}$ are all non-zero,
 therefore the linear functional $\int$ does not vanish identically on $M$, for
 example
$$\int((E^{ii}A \bar E^{ii})_\Lambda^J 
(\bar E^{ii}_\lambda B E^{ii}_\lambda)_{K'}^{\Gamma'})=
\frac{p_\lambda}{k_\lambda} (\bar C^{\otimes n} \bar E^{ii})_\Lambda^\Gamma 
E_\lambda^{ii}{}_K^J$$
is non-zero for a suitable choice of indices $K,J,\Gamma,\Lambda$. \hfill  $\Box$

\begin{thm} \label{thm1} Let $R$ and $\bar R$ be Hecke symmetries defined respectively
on $V$ and $\bar V$. Then there is a monoidal equivalence between 
$H$-comod and $\bar H$-comod sending $V$ to $\bar V$ and presvering
the braiding if and only if $R$ and $\bar R$ are defined with the same
parameter $q$ and have the same bi-rank.\end{thm}
\proof Assume that $R$ and $\bar R$ satisfies the condition of the theorem.
According to the lemma above it remains to prove that the monoidal functor
given by co-tensoring with $M$ send $V$ to $\bar V$ and $R$ to $\bar R$.
Indeed, by the definition of $V\Box_HM$,
 the map $\bar V \longrightarrow V \Box_H M$ given by
$\bar x_\lambda \longmapsto x_j \otimes a^j_\lambda$ is an 
injective $\bar H$-comodule homomorphism. According to Lemma \ref{lem3} and 
Schauenburg's result, 
$V \Box_H M$ is a simple $\bar H$-comodule, therefore
$\bar V$ is isomorphic to $V \Box_H M$. It is then easy to see that
$R$ is mapped to $\bar R$.

The converse statement is obvious. First, since $R$ is mapped to $\bar R$
they should be defined for the same value of $q$. Further, according
to Subsection \ref{sect1.2.4}, let $(r,s)$ and $(\bar r,\bar s)$ 
be the bi-ranks of $R$ and $\bar R$,
respectively. Then $\Gamma_n^{r,s}=\Gamma_n^{\bar r,\bar s}$ for all $n$,
whence $(r,s)=(\bar r, \bar s)$. \hfill  $\Box$

Notice that if $(r,s)\not=(\bar r, \bar s)$ and $r-s=\bar r-\bar s$ then 
$\Omega_n^{r,s}\cap\Omega_n^{\bar r,\bar s}= \emptyset $. This implies also that
the linear functional in Lemma \ref{lem3} is zero.

Theorem \ref{thm1} states that the study of comodules over a Hopf algebra associated
to a Hecke symmetry of bi-rank $(r,s)$ can be reduced to the study of the
Hopf algebra associated to the standard solution $R^{r,s}_q$. For the latter
Hopf algebra simple comodules were classified by R.B.~Zhang \cite{zhang}.
As an immediate consequence of Theorem \ref{thm1}, we have:
\begin{cor} Let $R$ be a Hecke symmetry of bi-rank $(r,s)$. Then the
homology of the associated Koszul complex (cf. Subsection  \ref{sect1.3.2})
 is concentrated at the term $K^{r,s}$ and has dimension one.
 Thus one has a homological determinant. \end{cor}
\proof In fact, the statement for $\bar R=R^{r,s}_q$ was proved by Manin \cite{manin2}.
Now, according to Theorem \ref{thm1}, for $M=M_{R,\bar R}$ the functor $-\Box_R M$
is fully faithful and exact hence the homology of the Koszul complex associated
to $R$ is concentrated at the term $r,s$ as the one associated to $\bar R$ is.
Since the homology group of the complex associated to $\bar R$ is one dimensional
and being an $\bar H$-comodule, it is an invertible comodule. Therefore
the homology group of the complex associated to $R$ is also invertible
as an $H$-comodule, hence is one-dimensional. \hfill $\Box$

\end{document}